\documentclass{article}
\usepackage{graphicx}

\def\be{\begin{equation}}
\def\ee{\end{equation}}
\def\ff#1{\mbox{\boldmath $#1$} }

\def\lam{\lambda}
\def\x{\ff{x}}

\def\vcode#1#2#3#4{\begin{figure}
\begin{center}
\begin{minipage}[c]{#1\textwidth}
{{\small #2 \hrule \vspace{5pt}   %
{\it #3}  \vspace{5pt} \hrule }}
\end{minipage}
\caption{#4}
\end{center}   \end{figure}} 

\begin{document}

\title{Flower Pollination Algorithm: A Novel Approach for Multiobjective Optimization}

\author{Xin-She Yang$^{a}$, Mehmet Karamanoglu$^{a}$, and Xingshi He$^{b}$  \\
${}^{a}$School of Science and Technology, \\
Middlesex University, The Burroughs, London NW4 4BT, UK \\
${}^{b}$College of Science, Xi'an Polytechnic University,  Xi'an, P. R. China
}

\date{}

\maketitle

\begin{abstract}
Multiobjective design optimization problems require multiobjective optimization
techniques to solve, and it is often very challenging to obtain high-quality
Pareto fronts accurately. In this paper, the recently developed flower pollination algorithm (FPA) is extended
to solve multiobjective optimization problems. The proposed method is used to solve
a set of multobjective test functions and two bi-objective design benchmarks, and
a comparison of the proposed algorithm with other algorithms has been made, which shows
that FPA is efficient with a good convergence rate. Finally, the importance for
further parametric studies and theoretical analysis are highlighted and discussed.
\end{abstract}

{\bf Citation Details:} 
X. S. Yang, M. Karamanoglu, X. S. He, Flower Pollination Algorithm: A Novel Approach for Multiobjective Optimization， {\it Engineering Optimization}, vol. 46, Issue 9, pp. 1222-1237 (2014).

\section{Introduction}

Real-world design problems in engineering and industry are usually multiobjective or multicriteria,
and these multiple objectives are often conflicting one another, which makes it impossible
to use any single design option without compromise. Common approaches are to provide
good approximations to the true Pareto fronts of the problem of interest
so that decision-makers can rank different options, depending on their preferences or
their utilities (Abbass and Sarker 2002; Babu and Gujarathi 2007; Cagnina et al. 2008;
Deb, 1999, 2000, 2001; Reyes-Sierra and Coello 2006).
Compared with single objective optimization, multiobjective optimization has its additional
challenging issues such as time complexity, inhomogeneity and dimensionality.
It is usually more time consuming to obtain the true Pareto fronts because
it usually requires to produce many points on the Pareto front for good approximations.

In addition, even accurate solutions on a Pareto front can be obtained, there is still
no guarantee that these solution points will distribute uniformly on the front.
In fact, it is often difficult to obtain the whole front
without any part missing.  For single objective optimization, the optimal solution can often
be a single point in the solution space, while for bi-objective optimization, the Pareto front
forms a curve, and for tri-objective cases, it becomes a surface.
In fact, higher dimensional problems can have
extremely complex hypersurface as its Pareto front (Madavan 2002; Marler and Arora 2004;
Yang 2010a; Yang and Gandomi 2012). Consequently, it is typically more challenging to solve
such high-dimensional problems.

In the current literature of engineering optimization,
a class of nature-inspired algorithms have shown their promising performance and
have thus become popular and widely used, and these algorithms are mostly swarm intelligence based
(Coello 1999; Deb et al. 2002; Geem et al. 2001; Geem 2009; Ray and Liew 2002;
Yang 2010,2010b,2011a; Gandomi and Yang 2011; Gandomi et al. 2012).
Metaheuristic algorithms such as particle swarm optimization, harmony search and cuckoo search are among the
most popular (Geem 2009; Yang 2010). For example, harmony search, developed by Zong Woo Geem in 2001
(Geem et al. 2001; Geem 2006, 2009), has been applied in many areas such as highly challenging water
distribution networks (Geem 2006) and discrete structural optimization (Lee et al. 2005).
Other algorithms such as shuffled frog-leaping algorithm and particle swarm optimizers have been
applied to various optimization problems (Eusuff et al. 2006; He et al. 2004; Huang 1996).
There are many reasons for the popularity of metaheuristic algorithms, and
flexibility and simplicity of these algorithms certainly contribute to their success.

The main aim of this paper is to extend the flower pollination algorithm (FPA),
developed by Xin-She Yang in 2012 (Yang 2012),
for single objective optimization to solve multiobjective optimization, and thus developed
a multi-objective flower pollination algorithm (MOFPA).
The rest of this paper is organized as follows: Section 2 outlines the basic characteristics of
flower pollination in nature and then introduce in detail the ideas of flower pollination algorithm.
Section 3 then presents the validation of the FPA by numerical experiments and a few selected multiobjective benchmarks.
Then, in Section 4, two real-world design benchmarks are solved to design a welded beam and a
disc brake, each with two objectives. Finally, some relevant issues are discussed  and conclusions are drawn in Section 5.

\section{Flower Pollination Algorithm}

\subsection{Characteristics of Flower Pollination}

It is estimated that there are over a quarter of a million types of flowering plants in Nature
and that about 80\% of all plant species are flowering species. It still remains
a mystery how flowering plants came to dominate the landscape from the Cretaceous period (Walker 2009).
Flowering plants have been evolving for at least more than 125 million years
and flowers have become so influential in evolution, it is unimaginable
what the plant world would look like without flowers. The main purpose of
a flower is ultimately reproduction via pollination. Flower pollination is
typically associated with the transfer of pollen, and such transfer
is often linked with pollinators such as insects, birds, bats and other animals.
In fact, some flowers and insects have co-evolved into a very specialized
flower-pollinator partnership. For example, some flowers can only attract and
can only depend on a specific species of insects or birds for successful pollination.

Pollination can take two major forms: abiotic and biotic. About 90\%
of flowering plants belong to biotic pollination. That is, pollen is
transferred by pollinators such as insects and animals. About 10\%
of pollination takes abiotic form which does not require any pollinators.
Wind and diffusion help pollination of such flowering plants,
and grass is a good example of abiotic pollination (ScienceDaily 2001; Glover 2007).
Pollinators, or sometimes called pollen vectors, can be very diverse.
It is estimated there are at least about 200,000 varieties of pollinators such as
insects, bats and birds. Honeybees are a good example of pollinators, and they
have also developed the so-called flower constancy. That is,
these pollinators tend to visit exclusive certain flower species while bypassing
other flower species. Such flower constancy may have evolutionary advantages
because this will maximize the transfer of flower pollen to the same or conspecific plants,
and thus maximizing the reproduction of the same flower species.
Such flower constancy may be advantageous for pollinators as well, because they
can be sure that nectar supply is available with their limited memory and minimum
cost of learning, switching or exploring. Rather than focusing on some unpredictable but potentially
more rewarding new flower species, flower constancy may require minimum investment
cost and more likely guaranteed intake of nectar (Waser 1986).

Pollination can be achieved by self-pollination or cross-pollination.
Cross-pollination, or allogamy, means pollination can occur from pollen
of a flower of a different plant, while self-pollination is the fertilization
of one flower, such as peach flowers, from pollen of the same flower or different flowers of the same plant,
which often occurs when there is no reliable pollinator available.
Biotic, cross-pollination may occur at long distance, and the pollinators such as
bees, bats, birds and flies can fly a long distance, thus they can considered as the
global pollination. In addition, bees and birds may behave as L\'evy flight behaviour
with jump or fly distance steps obeying a L\'evy distribution (Pavlyukevich 2007).
Furthermore, flower constancy can be considered as an increment step using the
similarity or difference of two flowers.

From the biological evolution point of view, the objective of the flower pollination
is the survival of the fittest and the optimal reproduction of plants in terms of
numbers as well as the most fittest. This can be considered as an optimization process of plant species.
All the above factors and processes of flower pollination interact so as to achieve
optimal reproduction of the flowering plants. Therefore, this may motivate us to design
new optimization algorithms.

\subsection{Flower Pollination Algorithm}

Flower pollination algorithm (FPA) was developed by Xin-She Yang in 2012 (Yang 2012),
inspired by the flow pollination process of flowering plants. FPA has been extended
to multi-objective optimization (Yang et al. 2013). For simplicity, the
following four rules are used:

\begin{enumerate}

\item Biotic and cross-pollination can be considered as a process of
global pollination, and pollen-carrying pollinators move in a way which
obeys L\'evy flights (Rule 1).

\item For local pollination, abiotic pollination and self-pollination are used (Rule 2).

\item Pollinators such as insects can develop flower constancy, which is equivalent to
a reproduction probability that is proportional to the similarity of two flowers involved (Rule 3).

\item The interaction or switching of local pollination and global pollination can
be controlled by a switch probability $p \in [0,1]$, slightly biased towards local
pollination (Rule 4).

\end{enumerate}

In order to formulate the updating formulas, the above rules
have to be converted into proper updating equations. For example, in the global pollination step,
flower pollen gametes are carried by pollinators such as insects,
and pollen can travel over a long distance because insects can often fly and move in a much
longer range. Therefore, Rule 1 and flower constancy (Rule 3)
can be  represented mathematically as
\be \x_i^{t+1}=\x_i^t + \gamma L(\lam) (\ff{g}_*-\x_i^t), \ee
where $\x_i^t$ is the pollen $i$ or solution vector $\x_i$ at iteration $t$, and $\ff{g}_*$ is
the current best solution found among all solutions at the current generation/iteration.
Here $\gamma$ is a scaling factor to control the step size.

\begin{figure}
\centerline{\includegraphics[height=3in,width=4in]{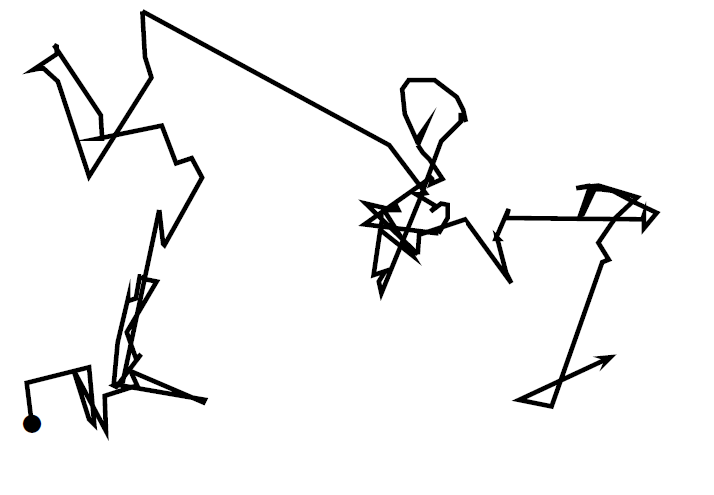}}
\caption{A series of 50 consecutive steps of L\'evy flights. \label{Levy-flight} }

\end{figure}

Here $L(\lam)$ is the parameter, more specifically the L\'evy-flights-based step size,
that corresponds to the strength of the pollination.
Since insects may move over a long distance with various distance steps, a L\'evy flight can be used
to mimic this characteristic efficiently. That is, $L>0$ is drawn from a L\'evy distribution
 \be L \sim \frac{\lam \Gamma(\lam) \sin (\pi \lam/2)}{\pi} \frac{1}{s^{1+\lam}}, \quad (s \gg s_0>0). \label{Levy-power} \ee
Here $\Gamma(\lam$) is the standard gamma function, and this distribution is valid for
large steps $s>0$. In theory, it is required that $|s_0| \gg 0$, but in practice
$s_0$ can be as small as $0.1$.  However, it is not trivial to generate pseudo-random step sizes that
correctly obey this L\'evy distribution (\ref{Levy-power}). There are a few methods for drawing
such random numbers, and the most efficient one from our studies is that the so-called
Mantegna algorithm for drawing step size $s$ by using two Gaussian distributions $U$ and $V$
by the following transformation (Mantegna 1994)
\be s=\frac{U}{|V|^{1/\lam}}, \quad U \sim N(0, \sigma^2), \quad V \sim N(0,1). \ee
Here $U \sim (0, \sigma^2)$ means that the samples are drawn from a Gaussian normal distribution
with a zero mean and a variance of $\sigma^2$. The variance can be calculated by
\be \sigma^2 =\Big[ \frac{\Gamma(1+\lam)}{\lam \Gamma((1+\lam)/2)} \cdot \frac{\sin (\pi \lam/2)}{2^{(\lam-1)/2}} \Big]^{1/\lam}. \ee
This formula looks complicated, but it is just a constant for a given $\lam$.
For example, when $\lam=1$, the gamma functions become $\Gamma(1+\lam)=1, \Gamma((1+\lam)/2)=1$ and
\be \sigma^2=\Big[\frac{1}{1 \times 1 } \cdot \frac{\sin (\pi \times 1/2)}{2^0}\Big]^{1/1}=1. \ee
It has been proved mathematically that the Mantegna algorithm can produce
the random samples that obey the required distribution (\ref{Levy-flight}) correctly (Mantegna 1994).
By using this pseudo-random number algorithm,  $50$ step sizes have been drawn to form
a consecutive 50 steps of L\'evy flights as shown in Fig. \ref{Levy-flight}.

\vcode{0.9}{{\sf Flower Pollination Algorithm (or simply Flower Algorithm) }} {
\indent Objective $\min$ or $\max f(\x)$, $\x=(x_1,x_2,..., x_d)$ \\
\indent Initialize a population of $n$ flowers/pollen gametes with random solutions  \\
\indent Find the best solution $\ff{g}_*$ in the initial population \\
\indent Define a switch probability $p \in [0,1]$ \\
\indent  {\bf while} ($t<$MaxGeneration) \\
\indent \qquad {\bf for} $i=1:n$ (all $n$ flowers in the population) \\
\indent \quad \qquad {\bf if} rand $<p$, \\
\indent \qquad \qquad Draw a ($d$-dimensional) step vector $L$ which obeys a L\'evy distribution \\
\indent \qquad \qquad Global pollination via $\x_i^{t+1}=\x_i^t + \gamma L (\ff{g}_*-\x_i^t)$ \\
\indent \qquad \quad {\bf else} \\
\indent \qquad \qquad Draw $\epsilon$ from a uniform distribution in [0,1] \\
\indent \qquad  \qquad Do local pollination via $\x_i^{t+1}=\x_i^t+\epsilon (\x_j^t-\x_k^t)$ \\
\indent \qquad \quad {\bf end if } \\
\indent \qquad \quad Evaluate new solutions \\
\indent \qquad \quad If new solutions are better, update them in the population \\
\indent \qquad {\bf end for} \\
\indent \qquad \quad Find the current best solution $\ff{g}_*$ \\
\indent  {\bf end while}  \\
\indent Output the best solution found }
{Pseudo code of the proposed Flower Pollination Algorithm (FPA). \label{fpa-code} }

For the local pollination, both Rule 2 and Rule 3 can be represented as
\be \x_i^{t+1}=\x_i^t + \epsilon (\x_j^t -\x_k^t), \ee
where $\x_j^t$ and $\x_k^t$ are pollen from different flowers of the same plant species.
This essentially mimics the flower constancy in a limited neighborhood. Mathematically,
if $\x_j^t$ and $\x_k^t$ comes from the same species
or selected from the same population, this equivalently becomes a local random walk if
$\epsilon$ is drawn from a uniform distribution in [0,1].

In principle, flower pollination activities can occur at all scales, both local and global.
But in reality,  adjacent flower patches or flowers in the not-so-far-away neighborhood
are more likely to be pollinated by local flower pollen than those far away.
In order to mimic this feature,  a
switch probability (Rule 4) or proximity probability $p$ can be effectively used to switch
between common global pollination to intensive local pollination.
To start with, a naive value of $p=0.5$ may be used as an initially value.
A preliminary parametric showed that $p=0.8$ may work better for most applications.

\subsection{Multiobjective Flower Pollination Algorithm (MOFPA)}

A multiobjective optimization problem with $m$ objectives can be written in general as
\be \textrm{Minimize }\;\; f_1(\x), f_2(\x), ..., f_m(\x),  \label{MOFPA-Obj} \ee
subject to the nonlinear equality and inequality constraints
\be h_j(\x)=0, \quad (j=1,2,...,J), \ee
\be g_k(\x) \le 0, \quad (k=1,2, ..., K). \ee
In order to use the techniques for single objective optimization or extend the methods for
solving multiobjective problems, there are different approaches to achieve this.
One of the simplest ways is to use a weighted sum to combine multiple objectives into a
composite single objective
\be f=\sum_{i=1}^m w_i f_i, \label{Single-Obj} \ee
with \be \sum_{i=1}^m w_i=1, \quad w_i > 0, \ee
where $m$ is the number of objectives and $w_i (i=1,...,m)$ are non-negative weights.

The fundamental idea of this weighted sum approach is that these weighting coefficients act
as the preferences for these multiobjectives. For a given set of $(w_1, w_2, ..., w_m)$,
the optimization process will produce a single point of the Pareto front of the problem.
For a different set of $w_i$, another point on the Pareto front can be generated.
With a sufficiently large number of combinations of weights, a good approximation to the
true Pareto front can be obtained. It is has proved that
the solutions to the problem with the combined objective (\ref{Single-Obj})
are Pareto-optimal if the weights are positive for all the objectives, and these
are also Pareto-optimal to the original problem (\ref{MOFPA-Obj}) (Miettinen 1999; Deb 2001).
In practice, a set of random numbers $u_i$ are first drawn for a uniform distribution $U(0,1)$.
Then, the weights $w_i$ can be calculated by normalization. That is
\be w_i=\frac{u_i}{\sum_{i=1}^m u_i}, \ee
so that $\sum_i w_i=1$ can be satisfied.
For example, for three objectives $f_1$, $f_2$ and $f_3$, three random numbers/weights can be drawn
from a uniform distribution $[0,1]$, and they may be $u_1=0.2915$, $u_2=0.9147$ and $u_3=0.6821$
in one instance of sampling runs. Then, we have $\sum_i=1.8883$, and $w_1=0.1544, w_2=0.4844, w_3=0.3612$. Indeed,
$\sum_i w_i=1.000$ is satisfied.

In order to obtain the Pareto front accurately with solutions relatively uniformly distributed
on the front, random weights $w_i$ should be used,  which should be as different as possible (Miettinen 1999).
From the benchmarks that have been tested, the weighted sum with random weights usually works well
as can be seen below.

\section{Validation and Numerical Experiments}

There are many different test functions for multiobjective optimization
(Zitzler and Thiele 1999; Ziztler et al. 2000; Zhang et al. 2009),
but a subset of some widely used functions provides a wide range of diverse properties
in terms of Pareto front and Pareto optimal set. To validate the proposed MOFPA,
a subset of these functions with convex, non-convex and discontinuous Pareto fronts
have been selected, including 7 single objective test functions and 4 multiobjective
test functions, and two bi-objective design problems.

\subsection{Single Objective Test Functions}
Before proceeding to solve multiobjective optimization problems,
the algorithm should first be validated by solving some well-known single objective test functions.
There are at least over a hundred well-known test functions. However, there is
no agreed set of test functions for validating new algorithms, though
some review and literature do exist (Ackley 1987; Floudas et al. 1999; Hedar 2012; Yang 2010a).
Here, a subset of seven test functions with diverse properties are used here.

The Ackley function can be written as
\be f_1(\x)=-20 \exp\Big[-\frac{1}{5} \sqrt{\frac{1}{d} \sum_{i=1}^d x_i^2}\Big] - \exp\Big[\frac{1}{d} \sum_{i=1}^d \cos (2 \pi x_i)\Big]
+20 +e, \ee
which has a global minimum $f_*=0$ at $(0,0,...,0)$.

The simplest of De Jong's functions is the so-called sphere function
\be f_2(\x) =\sum_{i=1}^n x_i^2, \quad
-5.12 \le x_i \le 5.12, \ee
whose global minimum is obviously $f_*=0$ at $(0,0,...,0)$. This function is unimodal and convex.

Easom's function
\be f_3(\x)=(-1)^{d+1} \prod_{i=1}^d \cos(x_i) \exp \Big[-\sum_{i=1}^d (x_i-\pi)^2  \Big], \ee
whose global minimum is $f_*=-1$ at $\x_*=(\pi, ..., \pi)$ within $-100 \le x_i \le 100$.
It has many local minima.

Griewank's function \be f_4(\x) = \frac{1}{4000} \sum_{i=1}^d x_i^2 - \prod_{i=1}^d \cos (\frac{x_i}{\sqrt{i}}) +1,
\quad -600 \le x_i \le 600, \ee
whose global minimum is $f_*=0$ at $\x_*=(0,0,...,0)$. This function is highly multimodal.

Rastrigin's function \be f_5(\x) = 10 d + \sum_{i=1}^d \Big[ x_i^2 - 10 \cos (2 \pi x_i) \Big],
\quad -5.12 \le x_i \le 5.12, \ee
whose global minimum is $f_*=0$ at $(0,0,...,0)$. This function is highly multimodal.

Rosenbrock's function \be f_6(\x) = \sum_{i=1}^{d-1} \Big[ (x_i-1)^2 + 100 (x_{i+1}-x_i^2)^2 \Big], \ee
whose global minimum $f_*=0$ occurs at $\x_*=(1,1,...,1)$ in the domain
$-5 \le x_i \le 5$ where $i=1,2,...,d$.

Zakharov's function
\be f_7(\x)=\sum_{i=1}^d x_i^2 +\Big(\frac{1}{2} \sum_{i=1}^d i x_i\Big)^2 +\Big(\frac{1}{2} \sum_{i=1}^d i x_i \Big)^4, \ee
has its global minimum $f_*=0$ at $(0,0,...,0)$.

In order to compare the performance of FPA with other existing algorithms,
each algorithm is first tested using the most widely used implementation and parameter settings.
For genetic algorithms (GA), a crossover rate of $p_{\rm crossover}=0.95$ and
a mutation rate of $p_{\rm mutation}=0.05$ are used (Holland 1975; Goldberg 1989; Yang 2010a).
For particle swarm optimization (PSO), a version with an inertia weight $\theta=0.7$ is used,
and its two learning parameters $\beta_1=\beta_2$ are set as $1.5$ (Kennedy and Eberhart 1995; Yang 2010a).
Also, to ensure a fair comparison,  the same population size should be used whenever possible.
So $n=25$ has been used for all three algorithms.

To get some insight into the parameter settings of the FPA,
a detailed parametric study has been carried out by varying $p$ from $0.05$ to $0.95$
with a step increase of $0.05$,  $\lam=1$, $1.25$, $1.5$, $1.75$, $1.9$ and $n=5, 10, 15, ..., 50$.
It has been found that $n=25$, $p=0.8$, $\gamma=0.1$ and $\lam=1.5$ works for most cases.
The parameter values used for all three algorithms are summarized in Table \ref{table-fpa-para}.

\begin{figure}
\centerline{\includegraphics[width=3in,height=2.5in]{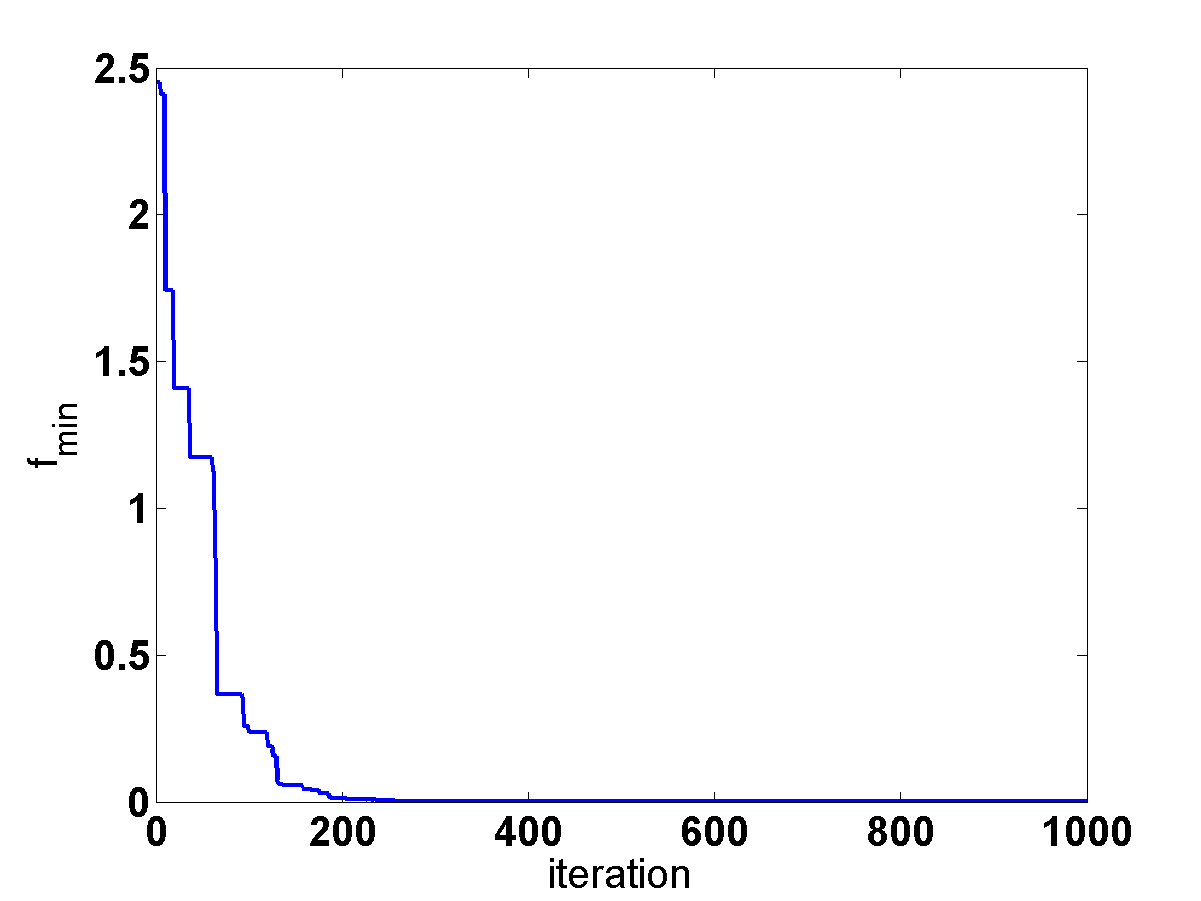}
\includegraphics[width=3in,height=2.5in]{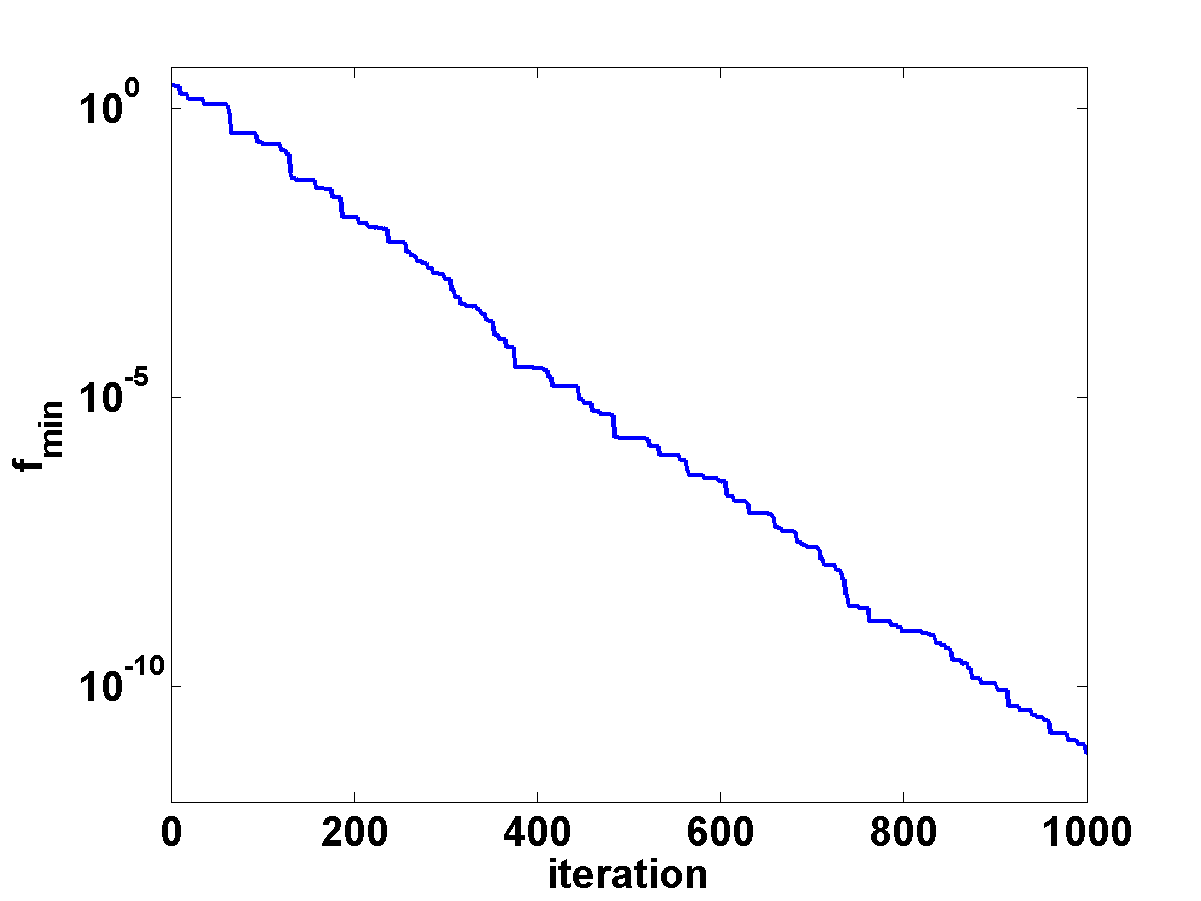} }
\caption{Convergence rate during iterations. The objective is plotted versus the iteration
(left), and the same results are shown in a logarithmic scale (right). \label{pafig-20} }
\end{figure}

\begin{table}[ht]

\caption{Parameter values for each algorithm. \label{table-fpa-para}}
\centering

\begin{tabular}{|l|l|} \hline
PSO & $n=25, \theta=0.7$, $\beta_1=\beta_2=1.5$ \\
\hline
GA & $n=25$, $p_{\rm crossover}=0.95, p_{\rm mutation}=0.05$ \\
\hline
FPA & $n=25$, $\lambda=1.5$, $\gamma=0.1$,  $p=0.8$  \\
\hline
\end{tabular}
\end{table}
The convergence behaviour of genetic algorithms and PSO during iterations
have been well studied in the literature. For FPA,
various statistical measures can be obtained from a set of runs.
For example, for the Ackley function $f_1$, the best objective values
obtained during each iteration can be plotted in a simple graph as
shown in Fig. \ref{pafig-20} where a logarithmic plot shows that
the convergence rate is almost exponential, which implies that
the proposed algorithm is very efficient.

For a fixed population size $n=25$, this is equal to  the total number of function evaluations is 25,000.
The best results obtained in terms of the means of the minimum values found are summarized
in Table \ref{table-fpa}.

\begin{table}[ht]

\caption{Comparison of algorithm performance: mean values. \label{table-fpa}}
\centering
\begin{tabular}{|r|r|r|r|}
\hline \hline
Functions & GA & PSO & FPA \\
\hline
$f_1$  & 8.29e-9  & 7.12e-12  & 5.09e-12 \\

$f_2$ & 6.61e-15 &  1.18e-24 & 2.47e-26 \\

$f_3$ & -0.9989 & -0.9998 & -1.0000 \\
$f_4$ & 5.72e-9 & 4.69e-9 & 1.37e-11 \\

$f_5$ & 2.93e-6 & 3.44e-6 & 4.52e-7 \\

$f_6$ & 8.97e-6 & 8.21e-8 & 6.19e-8 \\

$f_7$ & 8.77e-4 & 1.58e-4 & 9.53e-5 \\

\hline
\end{tabular}
\end{table}

\subsection{Multiobjective Test Functions}

In the rest of the paper, the parameters in MOFPA are fixed, based on
a preliminary parametric study, and $p=0.8$, $\lam=1.5$, and
a scaling factor $\gamma=0.1$ are used. The population size $n=50$ and the number of
iterations is set to $t=1000$. The following four functions will be tested:
\begin{itemize}
\item ZDT1 function with a convex front (Zitzler and Thiele 1999; Zitzler et al. 2000)
\[ f_1(x)=x_1, \quad f_2(x)=g (1-\sqrt{f_1/g}), \]
\be g=1+\frac{9 \sum_{i=2}^d x_i}{d-1}, \quad x_1 \in [0,1], \; i=2,...,30, \ee
where $d$ is the number of dimensions. The Pareto-optimality is reached when $g=1$.

\item ZDT2 function with a non-convex front
\[ f_1(x)=x_1, \quad f_2(x) =g (1-\frac{f_1}{g})^2, \]
where $g$ is the same as given in ZDT1.

\item ZDT3 function with a discontinuous front
\[ f_1(x) =x_1, \quad f_2(x)=g \Big[1-\sqrt{\frac{f_1}{g}}-\frac{f_1}{g} \sin (10 \pi f_1) \Big], \]
where $g$ in functions ZDT2 and ZDT3 is the same as in function ZDT1. In the ZDT3 function,
$f_1$ varies from $0$ to $0.852$ and $f_2$ from $-0.773$ to $1$.

\item LZ function (Li and Zhang, 2009; Zhang and Li, 2007)
\[ f_1=x_1 +\frac{2}{|J_1|} \sum_{j \in J_1} \Big [ x_j -\sin (6 \pi x_1 +\frac{j \pi}{d}) \Big]^2, \]
\be f_2=1-\sqrt{x_1} + +\frac{2}{|J_2|} \sum_{j \in J_2} \Big [ x_j -\sin (6 \pi x_1 +\frac{j \pi}{d}) \Big]^2, \ee
where $J_1=\{j|j$ is odd $\}$ and $J_2 =\{ j|j$ is even $\}$ where $2 \le j \le d$. This
function has a Pareto front $f_2=1-\sqrt{f_1}$ with a Pareto set
\be x_j=\sin (6 \pi x_1 + \frac{j \pi}{d}), \quad j=2,3, ..., d, \quad x_1 \in [0,1]. \ee

\end{itemize}

\begin{figure}[h]
\centerline{\includegraphics[height=2.7in,width=3in]{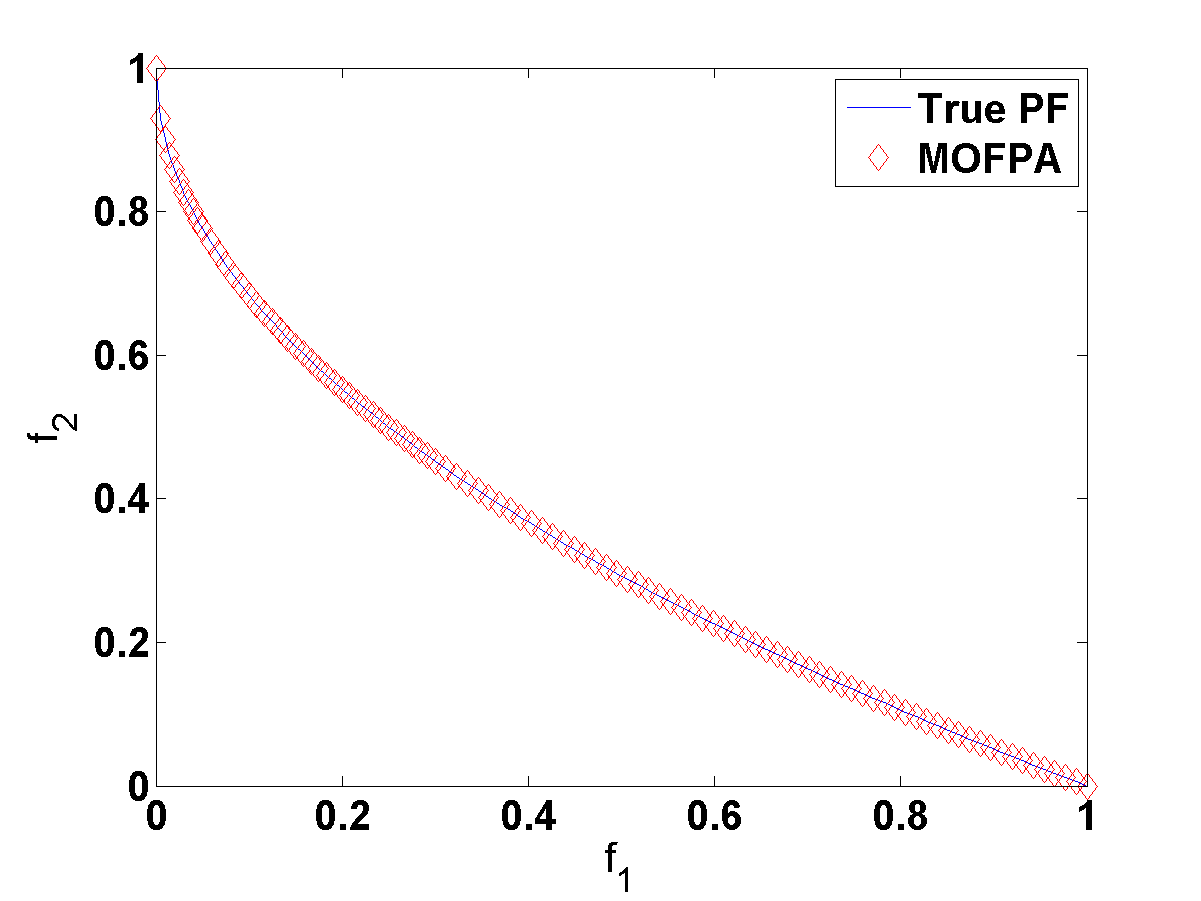}
\includegraphics[height=2.7in,width=3in]{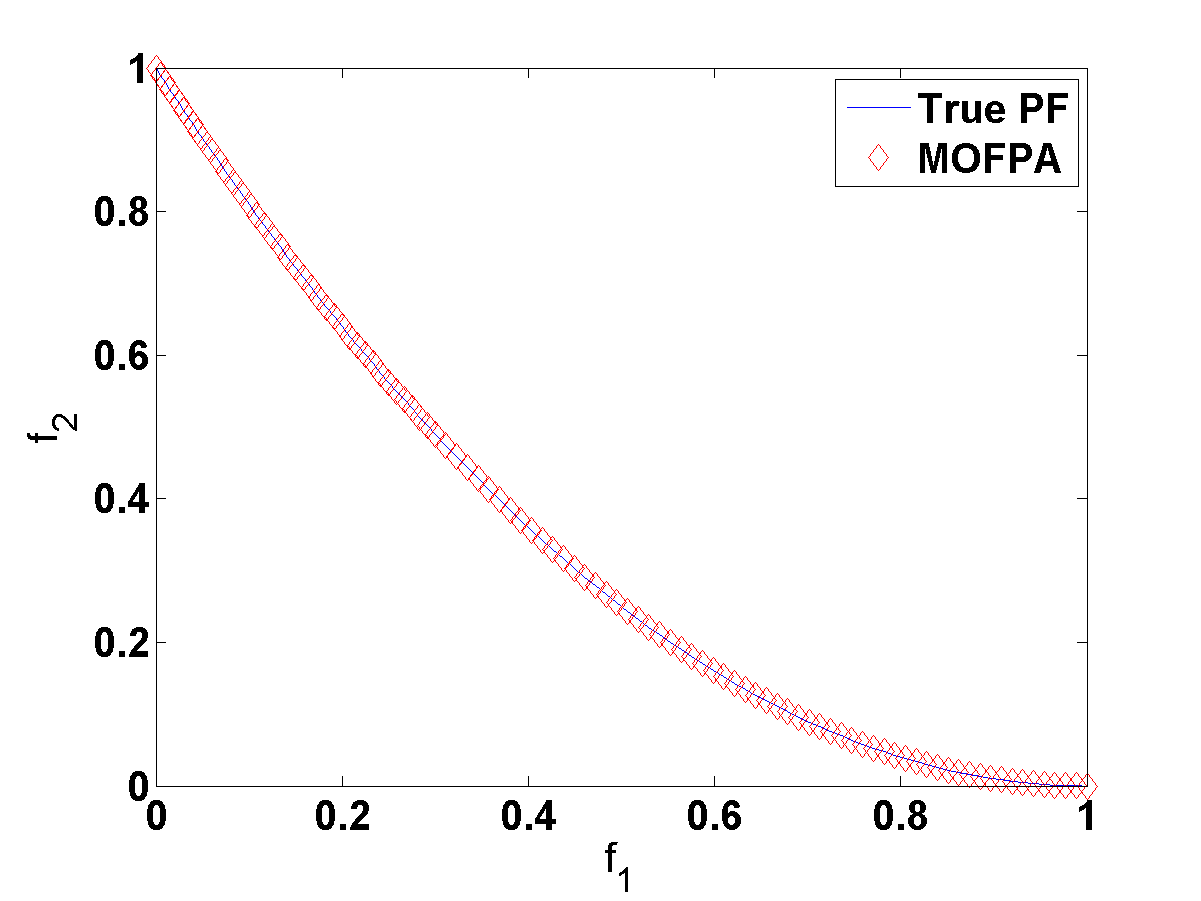} }
\caption{a) Pareto front of test function ZDT1, and b) Pareto front of test function ZDT2. \label{fig-300}}
\end{figure}

After generating 100 Pareto points by MOFPA, the Pareto front generated by MOFPA
is compared with the true front $f_2=1-\sqrt{f_1}$ of ZDT1 (see Fig. \ref{fig-300}).

Let us define the distance or error between the estimate Pareto front $PF^e$ to
its corresponding true front $PF^t$ as
\be E_f=||PF^e-PF^t||^2=\sum_{j=1}^N (PF_j^e-PF_j^t)^2, \ee
where $N$ is the number of points.

\begin{figure}[h]
\centerline{\includegraphics[height=2.7in,width=3in]{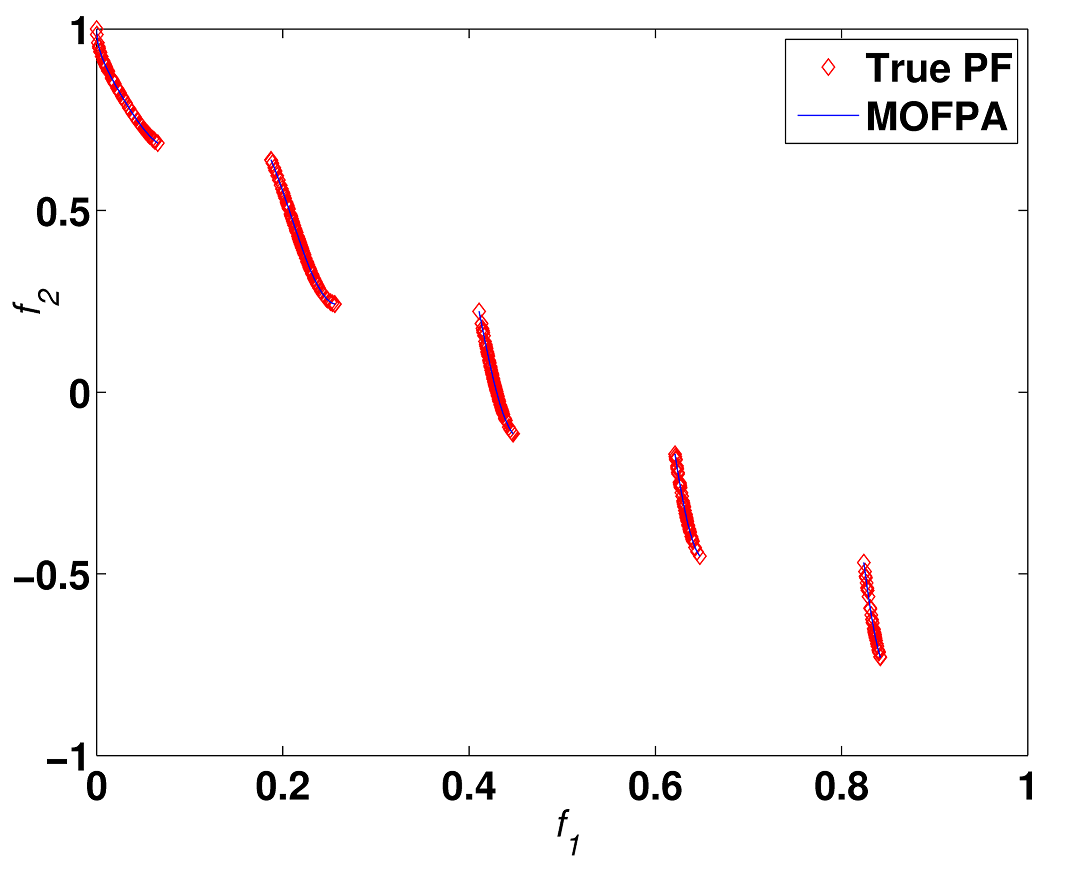}
\includegraphics[height=2.7in,width=3in]{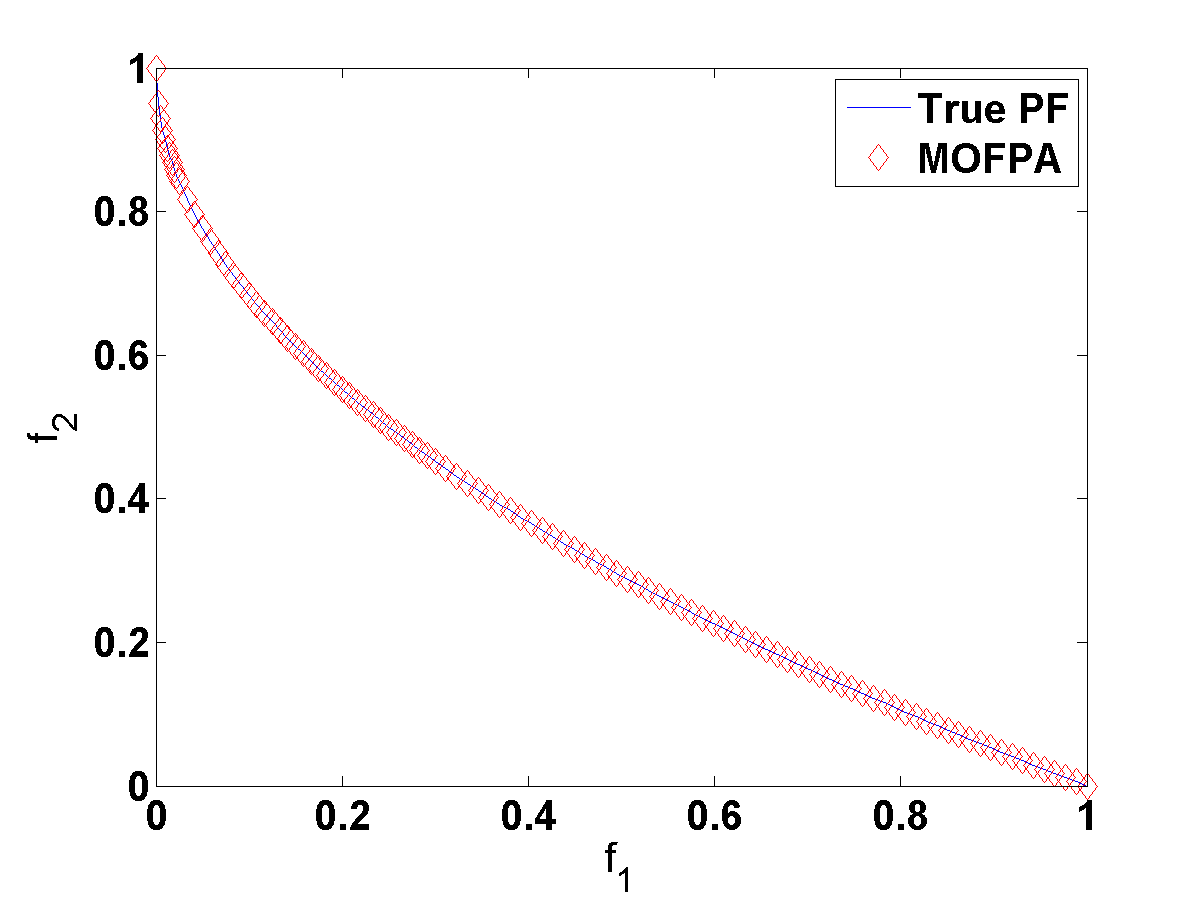} }
\caption{a) Pareto front of test function ZDT3, and b) Pareto front of test function LZ. \label{fig-400}}
\end{figure}

The variation of convergence rates or the convergence property can be viewed by
plotting out the errors during iterations. As this measure is an absolute measure,
which depends on the number of points. Sometimes, it is easier to use a
relative measure in terms of the generalized distance
\be D_g=\frac{1}{N} \sqrt{\sum_{j=1}^N (PF_j-PF_j^t)^2}. \ee
The results for all the functions are summarized in Table \ref{table-sum}, and the estimated
Pareto fronts and true fronts of other functions are shown in Fig. \ref{fig-300}
and Fig. \ref{fig-400}. In all these figures, the vertical axis is $f_2$ and
the horizontal axis is $f_1$.

\subsection{Analysis of Results and Comparison}

In order to compare the performance of the  proposed MOFPA with other established
multiobjective algorithms, we have selected a few algorithms with
available results from the literature. In case of the results are not available,
the algorithms have been implemented using well-documented studies
and then generated new results using these algorithms. In particular, other
methods are also used for comparison, including vector evaluated genetic algorithm
(VEGA) (Schaffer 1985), NSGA-II (Deb et al., 2000),
multiobjective differential evolution (MODE) (Babu 2007; Xue 2004),
differential evolution for multiobjective optimization (DEMO) (Robi\v c and Filipi\v c 2005),
multiobjective bees algorithms (Bees) (Pham and Ghanbarzadeh, 2007),
and strength pareto evolutionary algorithm (SPEA) (Deb et al. 2002; Madavan 2002).
The performance measures in terms of
generalized distance $D_g$ are summarized in Table \ref{table-sum-new} for all
the above major methods.

\begin{table}
\caption{Summary of results. \label{table-sum}}
\begin{center} \begin{tabular}{|l|l|l|}
\hline
Functions & Errors (1000 iterations) & Errors (2500 iterations) \\
\hline
ZDT1 & 1.1E-6   & 3.1E-19  \\
ZDT2 & 2.7E-6   & 4.4E-10  \\
ZDT3 & 1.4E-5   & 7.2E-12  \\
LZ & 1.2E-6  & 2.9E-12 \\
\hline
\end{tabular} \end{center} \end{table}

It is clearly seen from Table \ref{table-sum-new} that the proposed MOFPA obtained better
results for almost all four cases.

\begin{table}
\caption{Comparison of $D_g$ for $n=50$ and $t=500$ iterations. \label{table-sum-new}}
\begin{center} \begin{tabular}{|l|l|l|l|l|l|l|}
\hline
Methods &  ZDT1 & ZDT2 & ZDT3 & LZ \\
\hline \hline
VEGA & 3.79E-02 & 2.37E-03  & 3.29E-01 &   1.47E-03  \\
NSGA-II & 3.33E-02 & 7.24E-02  & 1.14E-01    & 2.77E-02  \\
MODE & 5.80E-03 & 5.50E-03 &     2.15E-02 & 3.19E-03  \\
DEMO & 1.08E-03 & 7.55E-04 &  1.18E-03 &    1.40E-03 \\
Bees  & 2.40E-02 & 1.69E-02 & 1.91E-01     & 1.88E-02 \\
SPEA & 1.78E-03 & 1.34E-03 & 4.75E-02   & 1.92E-03  \\ \hline
MOFPA & 7.11E-05 & 1.24E-05 & 5.49E-04   & 7.92E-05 \\
\hline
\end{tabular} \end{center} \end{table}

\section{Structural Design Examples}

Design optimization, especially design of structures, has many applications
in engineering and industry. As a result, there are many different benchmarks with
detailed studies in the literature (Kim et al. 1997; Pham and Ghanbarzadeh 2007; Ray and Liew 2002;
Rangaiah 2008).
In the rest of this paper, MOFPA will be used to solve two design case studies:
design of a beam and a disc brake (Osyczka and Kundu 1995; Ray and Liew 2002; Gong et al. 2009).

\subsection{Welded Beam Design}

Multiobjective design of a welded beam is a classical benchmark which has been solved by
many researchers (Deb 1999; Ray and Liew 2002). The problem has four design variables: the width $w$ and length $L$
of the welded area, the depth $d$ and thickness $h$ of the main beam. The objective is to minimize both
the overall fabrication cost and the end deflection $\delta$.

\begin{figure}
\centerline{\includegraphics[height=3in,width=4in]{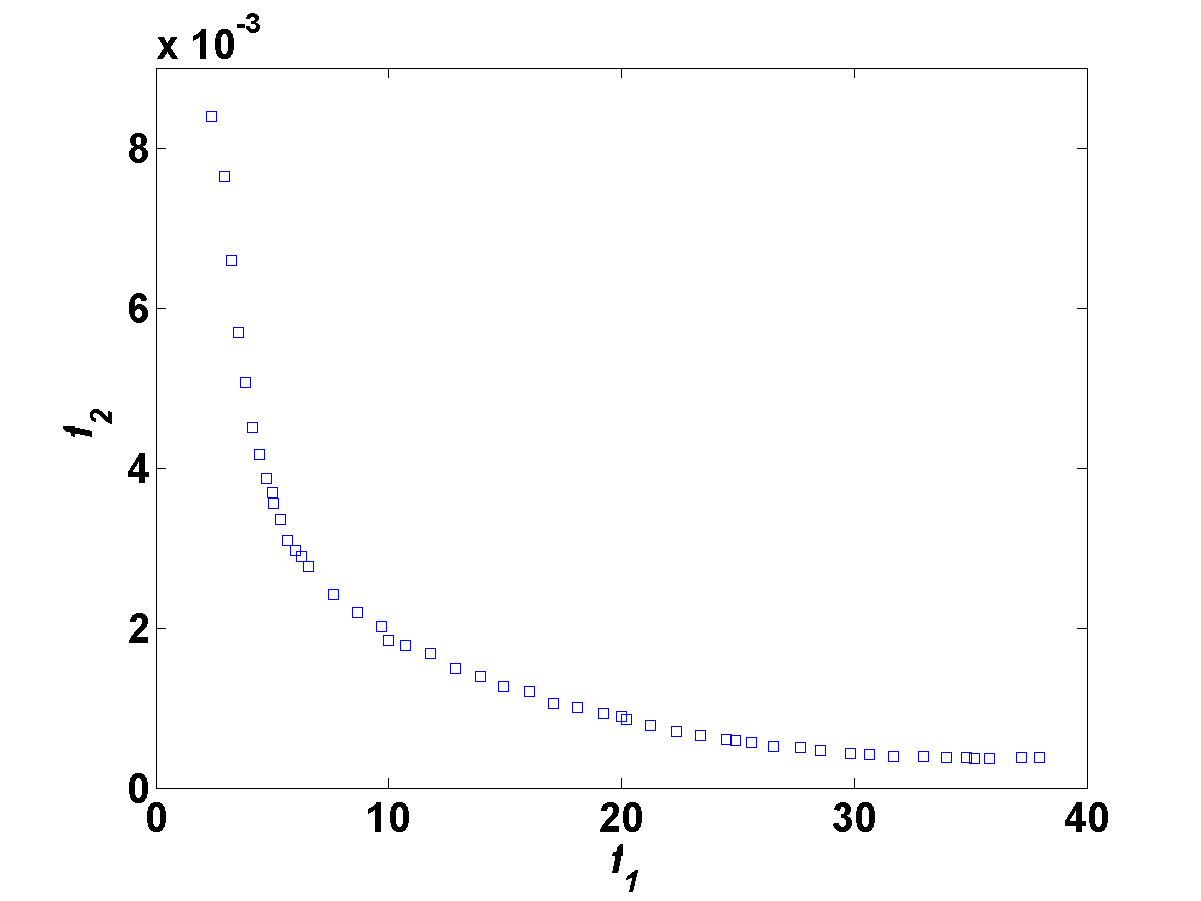} }
\caption{Pareto front for the bi-objective beam design where
the horizontal axis corresponds to cost and the vertical axis corresponds to deflection. \label{fig-440}}
\end{figure}

The detailed formulation can be found in (Deb 1999; Ray and Liew 2002; Gong et al. 2009).
Here  the main problem is rewritten  as
\be \textrm{minimise } \; f_1(\x)=1.10471 w^2 L + 0.04811 d h (14.0+L),
\; \textrm{ minimize } \; f_2=\delta, \ee
subject to
\be
\begin{array}{lll}
g_1(\x)=w -h \le 0, \vspace{2pt} \\ \vspace{3pt}
g_2(\x) =\delta(\x) - 0.25 \le 0, \\ \vspace{3pt}
g_3(\x)=\tau(\x)-13,600 \le 0, \\ \vspace{3pt}
g_4(\x)=\sigma(\x)-30,000 \le 0, \\ \vspace{3pt}
g_5(\x)=0.10471 w^2 +0.04811  h d (14+L) -5.0 \le 0, \\ \vspace{3pt}
g_6(\x)=0.125 - w \le 0, \\ \vspace{3pt}
g_7(\x)=6000 - P(\x) \le 0,
\end{array}
\ee
where
\be \begin{array}{ll}
 \sigma(\x)=\frac{504,000}{h d^2},  & Q=6000 (14+\frac{L}{2}), \\ \\
 D=\frac{1}{2} \sqrt{L^2 + (w+d)^2}, & J=\sqrt{2} \; w L [ \frac{L^2}{6} + \frac{(w+d)^2}{2}], \\ \\
 \delta=\frac{65,856}{30,000 h d^3}, &  \beta=\frac{QD}{J}, \\ \\
 \alpha=\frac{6000}{\sqrt{2} w L}, & \tau(\x)=\sqrt{\alpha^2 + \frac{\alpha \beta L}{D}+\beta^2}, \\ \\
 P=0.61423 \times 10^6 \; \frac{d h^3}{6} (1-\frac{d \sqrt{30/48}}{28}). &
\end{array} \ee
The simple limits or bounds are $0.1 \le L, d \le 10$ and
$0.125 \le w, h \le 2.0$.
This design problem has been solved using the MOFPA. The approximate Pareto front
generated by the 50 non-dominated solutions after 1000 iterations
are shown in Fig. \ref{fig-440}. This is
consistent with the results obtained by others (Ray and Liew 2002; Pham and Ghanbarzadeh 2007).

\subsection{Disc Brake Design}

The objectives are to minimize the overall mass and the braking time by choosing optimal design
variables: the inner radius $r$, outer radius $R$ of the discs, the engaging force $F$ and the number
of the friction surface $s$. This is under the design constraints such as the torque, pressure, temperature,
and length of the brake (Ray and Liew 2002; Pham and Ghanbarzadeh 2007).

This bi-objective design problem can be written as:
\be \textrm{Minimize } \; f_1(\x) =4.9 \times 10^{-5} (R^2-r^2) (s-1),
\quad f_2(\x)=\frac{9.82 \times 10^6 (R^2-r^2)}{F s (R^3-r^3)}, \ee
subject to
\be
\begin{array}{lll}
g_1(\x) =  20-(R-r) \le 0, \\[15pt]
g_2(\x) = 2.5 (s+1)-30 \le 0, \\[15pt]
g_3(\x) = \frac{F}{3.14 (R^2-r^2)} -0.4 \le 0, \\[15pt]
g_4(\x) = \frac{2.22 \times 10^{-3} F (R^3 -r^3)}{(R^2-r^2)^2} -1 \le 0, \\[15pt]
g_5(\x) =900- \frac{0.0266 F s (R^3-r^3)}{(R^2-r^2)} \le 0.
\end{array}
\ee
The simple limits are \be 55 \le r \le 80, \; 75 \le R \le 110, \; 1000 \le F \le 3000, \; 2 \le s \le 20. \ee
It is worth pointing out that $s$ is discrete. In general, MOFPA has to be extended in combination with
constraint handling techniques so as to
deal with mixed integer problems efficiently. However, since there is only one discrete variable,
the simplest branch-and-bound method is used here.

\begin{figure}[h]
\centerline{\includegraphics[height=4in,width=5.5in]{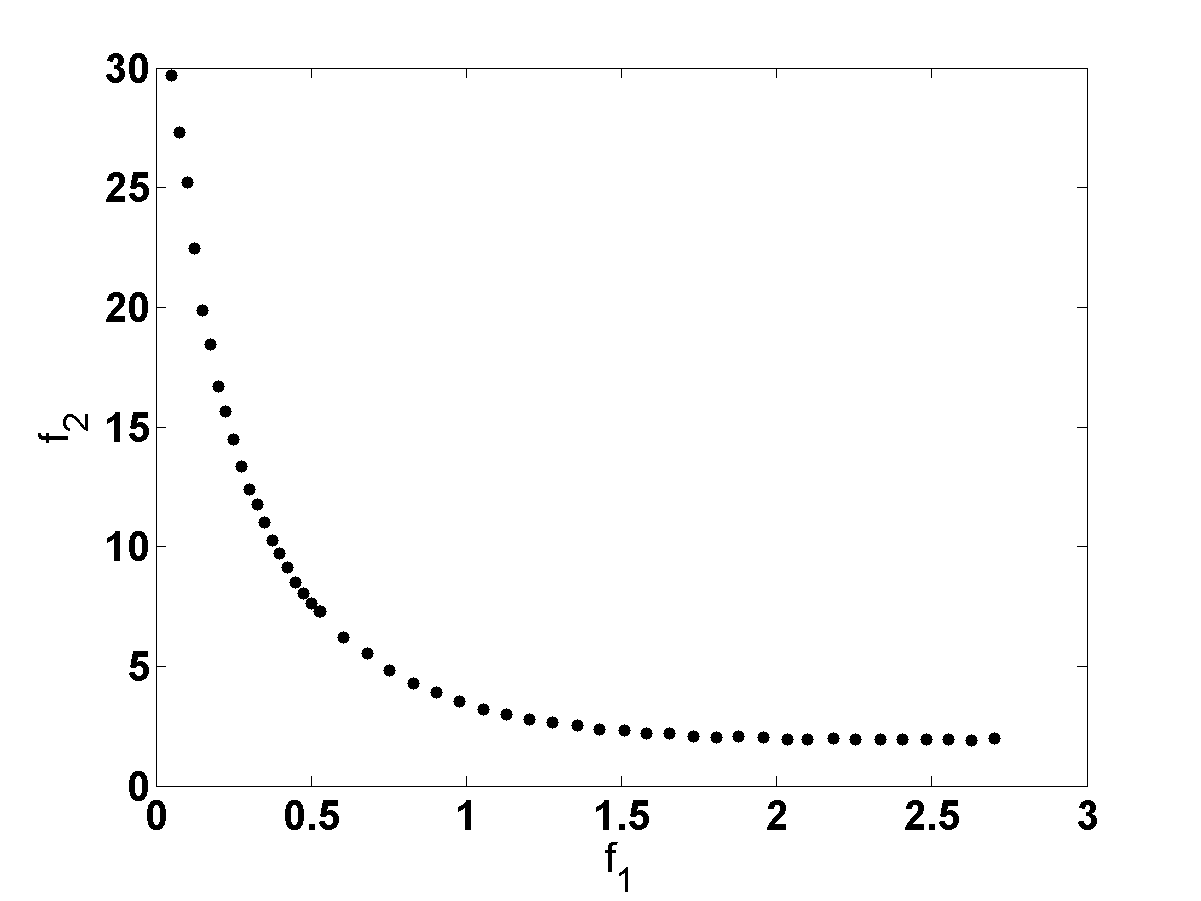} }
\caption{Pareto front of the disc brake design. \label{fig-800}}
\end{figure}

\begin{figure}[h]
\centerline{\includegraphics[height=3.5in,width=5in]{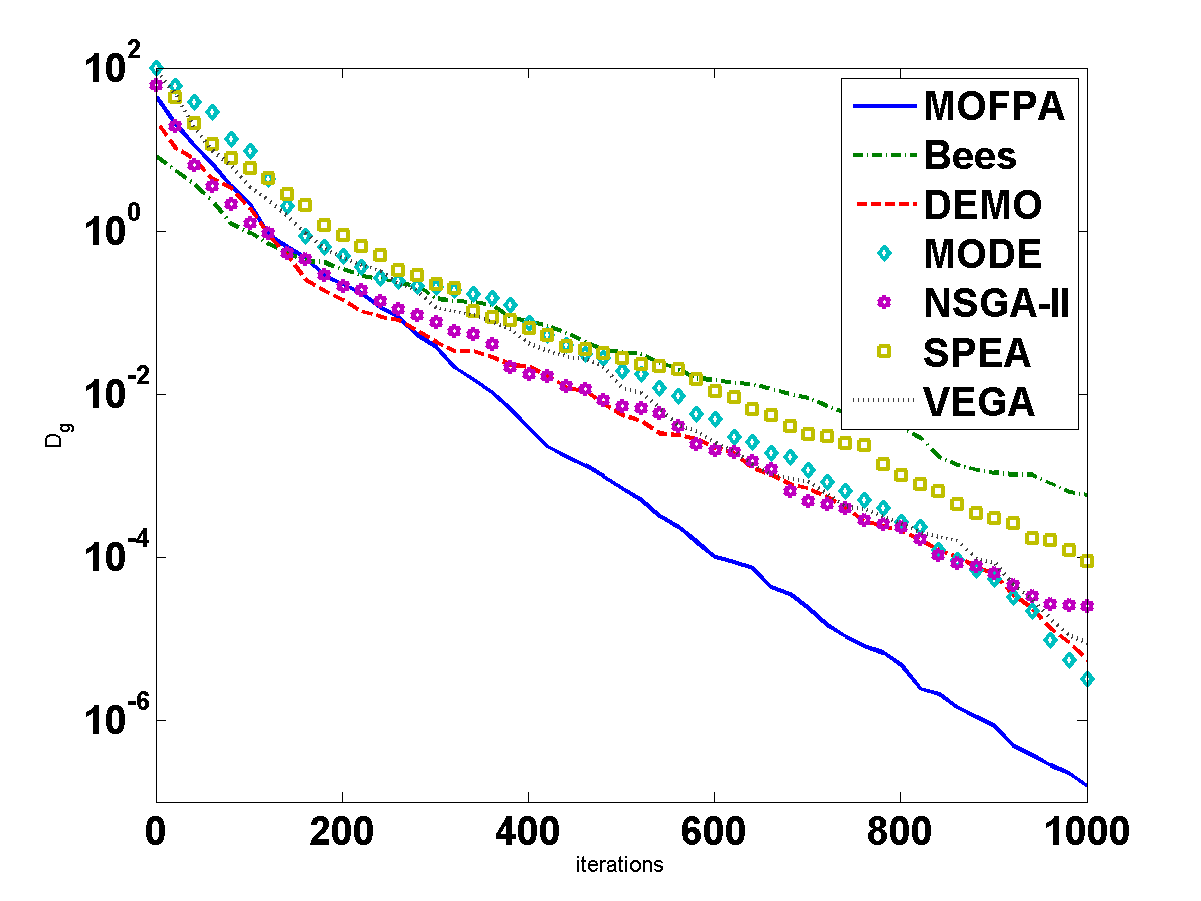} }
\caption{Convergence comparison for the disc brake design. \label{fig-900}}
\end{figure}

In order to see how the proposed MOFPA perform for the real-world design problems,
the same problem has also been solved using other available multiobjective algorithms.
50 solution points are geneated using MOFPA to form an approximate to the true Pareto front
after 1000 iterations, as shown in Fig. \ref{fig-800}.

The comparison of the convergence rates is plotted in the logarithmic scales
in Fig. \ref{fig-900}. It can be seen clearly that the convergence rate of MOFPA is the highest
in an exponentially decreasing way. This suggests that MOFPA provides better solutions in a more efficient way.

The above results for 11 test functions in total and two design examples
suggest that MOFPA is a very efficient algorithm for multiobjective optimization.
The proposed algorithm can deal with highly nonlinear, multiobjective optimization problems
with complex constraints and diverse Pareto optimal sets.

\section{Discussions and Conclusions}

Multiobjective optimization in engineering and industry is often
very challenging to solve, necessitating sophisticated techniques to tackle.
Metaheuristic approaches have shown promise and popularity in recent years.

In the present work, a new algorithm, called flower pollination algorithm, has been formulated
for multiobjective optimization applications by mimicking the pollination process of flowering plants.
Numerical experiments and design benchmarks have shown that the proposed algorithm
is very efficient with an almost exponential convergence rate, based on the
comparison of FPA with other algorithms for solving multiobjective optimization problems.

It is worth pointing out that mathematical analysis is highly needed in the future work
so as to gain insight into the true working mechanisms of the metaheuristic algorithms such as MOFPA.
FPA has the advantages such as simplicity and flexibility, and in many ways, it has
some similarity to that of cuckoo search and other algorithms with L\'evy flights (Yang 2010a, 2011b),
however, it is still not clear that why FPA works well. In terms of number of parameters,
FPA has only one key parameter $p$ together with a scaling factor $\gamma$,
which makes the algorithm easier to implement.
However, the nonlinearity in L\'evy flights make it difficult to analyse mathematically.
It can be expected that this nonlinearity in the algorithm formulations may be advantageous
to enhance the performance of an algorithm. More research may reveal the subtlety of this feature.

For multiobjective optimization, an important issue is how to ensure the solution points can
distribute relatively uniformly on the Pareto front for test functions. However, for real-world design
problems such as the design of a disc brake and a welded beam, the solutions are not quite uniform
on the Pareto fronts, and there is still room for improvement. However, simply generating more solution points
may not solve the Pareto uniformity problem easily. In fact, it is still a challenging problem on how to maintain a uniform
spread on the Pareto front, which requires more studies. It may be useful as a further research topic
to study other approaches for multiobjective optimization, such as the $\epsilon$-constraint method,
weighted metric methods, Benson's method, utility methods, and evolutionary methods (Miettinen 1999;
Coello 1999; Deb 2001).

On the other hand, further studies can focus on more detailed parametric analysis and gain insightful
of how algorithm-dependent parameters can affect the performance of an algorithm.
Furthermore, the linearity in the main updating formulas makes it possible to do some theoretical
analysis in terms of dynamic systems or Markov chain theories, while the nonlinearity in terms of L\'evy flights
can be difficult to analyze FPA exactly. All these can form useful topics for further research.


\end{document}